\documentclass[12pt]{amsart}

\usepackage{amsfonts,amsmath,latexsym,amssymb,verbatim,amsbsy}
\usepackage{amsthm}

\usepackage{pstricks}

\theoremstyle{plain}
\newtheorem{THEOREM}{Theorem}[section]

\newtheorem{theorem}[THEOREM]{Theorem}

\theoremstyle{definition}

\theoremstyle{remark}

\newcommand{\thm}[1]{Theorem~\ref{#1}}

\newcommand{\R}{\ensuremath{\mathbb{R}}}   

\def \a {\alpha}
\def \b {\beta}
\def \d {\delta}
\def \g {\gamma}
\def \e {\varepsilon}

\def \n {\nabla}
\def \s {\sigma}

\def \t {\tau}

\def \calD {\mathcal{D}}
\def \D {\Delta}

\def \backs {\backslash}

\def \< {\langle}
\def \> {\rangle}
\def \p {\partial}
\def \ra {\rightarrow}
\def \ss {\subset}

\newcommand{\be}{\begin{equation}}
\newcommand{\ee}{\end{equation}}

\def \sp {\s_{\text{p}}}
\def \loc {\mathrm{loc}}

\DeclareMathOperator{\supp}{supp} %
\DeclareMathOperator{\diver}{div} %
\DeclareMathOperator{\tr}{Tr} %

\begin{document}

\title[Energy equality for 3D NSE]{A geometric condition implying energy equality for solutions of 3D Navier-Stokes equation.}

\author{R. Shvydkoy}
\thanks{The work was partially supported by NSF grant DMS -- 0604050. The author is grateful to A. Cheskidov and G. Seregin for stimulating discussions.}
\address{Department of Mathematics, Stat. and Comp. Sci.\\
 851 S  Morgan St., M/C 249\\
        University of Illinois\\
        Chicago, IL 60607}
\email{shvydkoy@math.uic.edu}

\date{\today}

\begin{abstract}
We prove that every weak solution $u$ to the 3D Navier-Stokes
equation that belongs to the class $L^3L^{9/2}$ and $\n u$ belongs
to $L^3L^{9/5}$ localy away from a $1/2$-H\"{o}lder continuous
curve in time satisfies the generalized energy equality. In particular
every such solution is suitable.
\end{abstract}



\maketitle

\section{Introduction}

In this note we discuss the energy balance equality for weak
solutions of the 3D Navier-Stokes equations. The system of the
Navier-Stokes equations is given by
\begin{align}
 u_t + \diver(u\otimes u) + \n p & = \nu \D u ,\label{nse}\\
 \n \cdot u & = 0, \label{diver}
\end{align}
where $u$ is the velocity field, $p$ the internal pressure.
We focus primarily on the case of $\R^3$. A weak solution to \eqref{nse}--\eqref{diver}
is a pair of distributions $(u,p) \in \calD'((0,T)\times \R^3)^4$ with $u \in L^2((0,T)\times \R^3)_{\text{loc}}$ such that
\begin{align}
 -\iint u \psi_t - \nu \iint u  \D \psi &= \iint \left( \tr[(u \otimes u)\cdot \n \psi] +
p \diver{\psi}\right)\label{wnse} \\
\iint u \cdot \n \phi &= 0,
\end{align}
for all $(\psi,\phi) \in \calD((0,T) \times \R^3)^4$, where
$\calD$ stands for the space of $C^\infty$-smooth compactly
supported functions. The classical existence theorem of Leray
\cite{Leray} states that given any divergence free initial
condition $u_0 \in L^2(\R^3)$ one can find at least one weak
solution $(u,p)$ with $u \in L^\infty L^2 \cap L^2H^{1}$ (here and
throughout $L^rX = L^r([0,T];X(\R^3))$), $u(t) \ra u_0$ strongly
in $L^2$ as $t \ra 0$, the pressure is given by
\begin{equation}\label{press}
 p = \sum_{i,j = 1}^3R_iR_j(u_i u_j),
\end{equation}
where $R_i$'s are the classical Riesz transforms, and the following energy inequality
\begin{equation}\label{ei}
 \frac{1}{2} \int_{\R^3 \times \{t\}} |u|^2 + \nu \int_{t_0}^t \int_{\R^3} |\n u |^2 \leq \frac{1}{2} \int_{\R^3 \times \{t_0\}} |u|^2
\end{equation}
holds for all $t \in (0,T]$ and a.e. $t_0 \leq t$ including $t_0 = 0$.
Following Serrin \cite{Serrin} one can further tune $u$ on a time set of measure zero
to achieve weak continuity in $L^2$. We denote $u \in C_w([0,T]; L^2(\R^3)) = C_wL^2$.

The lack of exact equality in \eqref{ei} is a pitiful  deficiency
of Leray's solutions that to date remains unresolved. The main
difficulty arises in the fact that the mollified velocity field
$u_\d$ may have a non-vanishing energy flux due to the nonlinear
term, i.e. the equality
\begin{equation}\label{flux}
\lim_{\d \ra 0} \iint \tr[(u\otimes u)_\d \cdot \n u_\d] = 0
\end{equation}
may fail. In \cite{Lions} Lions showed that if $u \in L^4L^4$ then
the energy equality holds. Technique developed in \cite{LSU} for
general parabolic equations reproduces Lions' result as well. By
interpolation with $L^\infty L^2$ one automatically obtains the
range of conditions
\begin{equation}\label{shinbrot}
u \in L^r L^s, \text{ for some } 2/r + 2/s = 1,\; s \geq 4.
\end{equation}
These were shown to work in  any spacial dimension by Shinbrot
\cite{Shinbrot}.  The above results are based on proving
continuity of the trilinear form $u\cdot \n v \cdot w$ in
$L^\infty L^2 \cap L^2 H^1 \cap L^r L^s$ making it possible to
carry out the standard mollification argument in order to obtain
\eqref{flux}. A different approximation procedure was proposed by
Kukavica \cite{Kukavica}. It allowed to essentially use
cancellations present in the nonlinear term. The extra regularity
requirement was thus relieved from the velocity field and put to the pressure
giving the condition $p \in L^2((0,T)\times \R^3)_{\text{loc}}$.
In view of \eqref{press} this is a weaker yet dimensionally the
same version of Lions' condition.

After a recent progress on Onsager's conjecture for the Euler
equations (see \cite{ccfs,CET,Duchon}) dimensionally sharper
conditions were found in \cite{ccfs,cfs}. Namely, the energy
equality holds if $u \in L^3 B^{1/3}_{3,p}$ in the case of $\R^3$
or periodic domain, or if $u \in L^3 D(A^{5/12})$ in the case of a
bounded domain. Here $B^{1/3}_{3,p}$ is the Besov space with
smoothness $1/3$ integrability $3$ and summability $1\leq p <
\infty$, and $A$ denotes the Stokes operator. The dimensional
$L^p$-analogue of these spaces is $L^3L^{9/2}$, which lies outside
of \eqref{shinbrot}. In fact the cube of the dimension of
$L^3L^{9/2}$ is the same as the dimension of the energy flux
in \eqref{flux}, suggesting that this space might be the optimal one for
any argument based on direct control of the flux. We thus
conjecture that every weak solution to \eqref{nse}--\eqref{diver}
in the class $u \in C_wL^2 \cap L^2H^1 \cap L^3L^{9/2}$ verifies
the energy equality. In this paper we prove the following result
in this direction.

\begin{theorem}\label{T:main}
Let $s \in C^{1/2}([0,T]; \R^3)$ and $(u,p)$ be a weak solution to the NSE satisfying the following conditions
\begin{itemize}
 \item[(i)] $u\in C_w L^2 \cap L^2H^1 \cap L^3L^{9/2}$;
\item[(ii)] $\n u \in L^3L^{9/5}((0,T)\times \R^3 \backs \mathrm{Graph(s)})_{\loc}$,
\end{itemize}
and $p$ is given by \eqref{press}. Then $(u,p)$ satisfies the generalized energy equality:
\begin{multline}\label{mgee}
\int_{\R^3\times\{t\}} |u|^2 \phi + 2\nu \int_0^t\int_{\R^3} |\n
u|^2 \phi =\int_{\R^3\times\{0\}} |u|^2 \phi +\\ +
\int_0^t\int_{\R^3} \left[ |u|^2(\phi_t+ \nu \D \phi)+
(|u|^2+2p)u\cdot \n \phi\right],
\end{multline}
for all $\phi \in \calD([0,T]\times \R^3)$, and $t\in [0,T]$.
\end{theorem}
A few remarks are in order. First, we note that for any weak solution $(u,p)$ with $u \in L^\infty L^2 \cap L^2H^1$, the pressure is restored via
\eqref{press} up to an $x$-independent distribution. Thus, any condition $p \in L^rL^s$ would already imply \eqref{press}, and hence the inclusions
\begin{equation}\label{p-est}
 p \in L^rL^s, \text{ for } 2/r + 3/s = 3,\quad 1<s \leq 3.
\end{equation}
This justifies the integrals in \eqref{mgee}. Furthermore,
considering the sequence $\phi(x/R)$ with $\phi=1$ near $0$ and
$R\ra 0$ we recover the global energy equality.

Following  \cite{CKN,Lin} weak solutions with the  natural bounds
on $u$ and $p$ satisfying the generalized energy inequality in
\eqref{mgee} for all non-negative $\phi$ are called suitable. The
end result of the partial regularity theory developed in
\cite{CKN,Scheffer3,Scheffer,Scheffer2} showed that the set of
singular points of every such solution has zero one-dimensional
parabolic Hausdorff measure. From this point of view the
H\"{o}lder condition on $s$ seems rather generous as the graph of
a $C^{1/2}$-function may have parabolic dimension as large as $2$.
Yet it is essential for the argument that the curve $s$ is
extended in time. For instance, assuming that $s$ is a smooth
curve on a slice $\R^3 \times \{t_0\}$ our argument necessitates the stronger
condition $u\in L^3 L^6$, which already falls into the range of
\eqref{shinbrot}.

Finally let us note that condition $(ii)$ is dimensionally the
same as the $L^3L^{9/2}$ condition on $u$. As such it lies out of
reach of the classical Prodi-Serrin condition \cite{Serrin} or the condition
proposed in \cite{Veiga}. So, theoretically by requiring $(ii)$ we
do not exclude the possibility of having singularities away from
the graph of $s$. In Section \ref{S:ext} we will continue our
discussion of \thm{T:main}.

\section{Proof of \thm{T:main}}

 The proof is based on an approximation procedure.
So, let us fix a mollifier $\b \geq 0$, $\b \in C_0^\infty(B_1)$ with
$\int \b = 1$, where $B_\rho = \{|x| \leq \rho\}$. For a distribution $u\in \calD'(\R^3)^n$ we denote
$$
u_\d(x) = \d^{-3} \int_{\R^3} \b \left(y\d^{-1} \right)
u(x-y) dy.
$$
If $u$ is weak solution to the NSE and $u \in C_w([0,T];L^2)$, then we have
\begin{multline}\label{weakNSE2}
  \int_{\R^3 \times \{t\}} u \psi -\iint u \psi_t - \nu \iint u \D \psi \\= \int_{\R^3 \times \{0\}} u \psi +\iint \left( \tr[(u \otimes u)\cdot \n \psi] +
p \diver{\psi}\right)
\end{multline}
holds for all $t\in [0,T]$ and $\psi \in \calD([0,T]\times \R^3)^3$. Substituting  $$\psi = \d^{-3}\b((x-\cdot)\d^{-1}) e_k,$$
where $e_k$ are the vectors of the standard unit basis, we immediately obtain
$$
u_\d(t) = u_\d(0) + \left\{ \int_0^t[ u_k(s)\ast \b'+ (u_k(s)u_j(s))\ast \b''_j + p(s)\ast\b_k''']ds \right\}_{k=1}^3,
$$
for some $\b',\b''_j,\b'''_j \in \calD(\R^3)$. In view of $u\in
L^\infty L^2\cap L^2L^6$ and \eqref{p-est}, the function under the
integral belongs to $L^rX$, where $X$ is any Sobolev space
$W^{\g,p}$, $s\geq 0$, $p\geq 2$, and $r < \infty$. This implies
that $u_\d$ is absolutely continuous in $X$ with Freschet
derivative $\p_t u_\d \in L^rX$. By the standard approximation
argument, functions with such smoothness are allowed in
\eqref{weakNSE2}. We therefore can substitute a test-function of
the form $\psi = (u_\d \phi)_\d$, where $\phi \in
\calD([0,T]\times \R^3)$.

We now proceed with the construction of the appropriate test-function. In order to cut off
the graph of $s$ we first extend $s$ beyond $[0,T]$ by defining
\begin{equation}
s^{\mathrm{ext}}(t) = \left\{
  \begin{array}{ll}
    s(0), & t<0; \\
    s(t), & 0\leq t <T; \\
    s(T), & t\geq T.
  \end{array}
\right.
\end{equation}
Clearly, $s^{\mathrm{ext}} \in C^{1/2}(\R; \R^3)$. Second, we define
\begin{equation}
s_\e(t) = \int_{\R} \e^{-2} \a(\t \e^{-2}) s^{\mathrm{ext}}(t-\t) d
\t,
\end{equation}
for some mollifier $\a$ and $\e>0$. The following approximation inequalities easily follow:
\begin{align}
\sup_{0\leq t\leq T} |s_{\e}(t) - s(t)| &<\e;\label{appr1}\\
\sup_{0\leq t\leq T} |s'_{\e}(t)| \leq 1/\e.
\end{align}
Next we introduce a cut-off function $\chi \in C^\infty(\R^3)$
with $\chi \geq 0$, $\chi \equiv 0$ in $B_2$ and $\chi \equiv 1$
in $\R^3 \backs B_3$. Denote
$$
\chi_\e(x,t) = \chi\left(\frac{x-s_\e(t)}{\e}\right).
$$
Notice that in view of \eqref{appr1},
\begin{equation}\label{supp}
\supp{\chi_\e} \ss \{(x,t): |x-s(t)|>\e\},
\end{equation}
and $\chi_\e$ in infinitely smooth in time-space. Let us note the following inequality
\begin{equation}\label{der}
 \sup_{t\in[0,T]} \|D_x^\g \chi_\e\|_p \sim \e^{\frac{3}{p} - |\g|},
\end{equation}
for any $1\leq p\leq \infty$ and multiindex $\g$.

Let us fix an arbitrary $\phi \in \calD([0,T]\times \R^3)$ and define the following test-function
\begin{equation}
 \psi = (u_\d \phi \chi_\e)_\d.
\end{equation}
As we substitute this function into \eqref{weakNSE2} we will adhere
to the same order of limits as $\e,\d \ra 0$ in all our subsequent
computations. Namely, first $\d \ra 0$ and then $\e\ra 0$. Let us
assign letters to the terms of equation \eqref{weakNSE2} by writing it
as
\begin{equation}\label{assign}
A-B-C = D + E.
\end{equation}
We now examine each term separately.

First, let us notice that integration by parts carried out in $B$ results in appearance of two terms that cancel out with $A$ and $B$ plus the following
$$
\frac{1}{2}\int_{\R^3\times\{t\}} |u_\d|^2 \phi \chi_\e -
\frac{1}{2}\int_{\R^3\times\{0\}} |u_\d|^2 \phi \chi_\e -
\frac{1}{2}\int_0^t \int_{\R^3} |u_\d|^2 (\phi \chi_\e)_t.
$$
The first two integrals converge to the corresponding terms in
\eqref{mgee}, while the third integral is given by
$$
\frac{1}{2}\iint |u_\d|^2 (\phi \chi_\e)_t =
\frac{1}{2}\iint|u_\d|^2 \phi_t \chi_\e +
\frac{1}{2}\iint |u_\d|^2 \phi\ s'(t)\cdot \n
\chi_\e.
$$
Clearly, the first integral on the right hand side converges to its
natural limit
$$
\frac{1}{2}\iint |u|^2 \phi_t
$$
producing the corresponding term in \eqref{mgee}. The second
integral converges to zero. Indeed, by \eqref{der} and H\"{o}lder, one obtains
\begin{align*}
&\left|\iint |u_\d|^2 \phi\ s'(t)\cdot \n \chi_\e
\right| \leq \int_0^t \left(  \int_{|x-s_\e(t)| \leq 3\e}
|u_\d|^6 |\phi|^3 \ dx \right)^{1/3} |s'_\e(t)| \e\ dt \\
& \leq \int_0^t \left(  \int_{|x-s_\e(t)| \leq 3\e} |u_\d|^6 \
dx \right)^{1/3} dt \leq \int_0^t \left( \int_{|x-s_\e(t)|
\leq 3\e+\d } |u|^6 \ dx \right)^{1/3} dt,
\end{align*}
and the latter converges to zero as $\d,\e \ra 0$ since $u\in L^2L^6$.

Let us now examine term $C$. We have
\begin{align*}
-C &= \nu \iint \n u_\d \cdot \n( u_\d \phi \chi_\e) =
\nu \iint |\n u_\d|^2 \phi \chi_\e - \frac{\nu}{2}
\iint |u_\d|^2 \D(\phi \chi_\e) \\
& = \nu \iint |\n u_\d|^2 \phi \chi_\e -
\frac{\nu}{2} \iint |u_\d|^2 \D\phi \chi_\e \\
& - \nu \iint |u_\d|^2 \n\phi \cdot \n \chi_\e -
\frac{\nu}{2} \iint |u_\d|^2 \phi \D\chi_\e \\
& = C_1 - C_2 - C_3 - C_4.
\end{align*}
By the standard convergence theorems we see that $C_1$ and $C_2$
converge to the corresponding terms in \eqref{mgee}, while in view of \eqref{der}
\begin{align*}
|C_3| & \leq \e \int_0^t \|u_\d \n \phi\|_6^2 dt \ra 0, \\
|C_4| & \leq \int_0^t \left(  \int_{|x-s_\e(t)| \leq 3\e}
|u_\d|^6 |\phi|^3 \ dx \right)^{1/3} dt \ra 0.
\end{align*}

Let us examine term $E$. We have
$$
E = \iint \tr[ (u\otimes u)\cdot \n \psi ] +
\iint p\diver \psi = F+G.
$$
We can write
\begin{align}
 F& = \iint \tr[ (u\otimes u)_\d \cdot \n (u_\d \phi \chi_\e) ] \\
& = \iint \tr[ r_\d(u,u) \cdot \n (u_\d \phi \chi_\e) ] \\
& + \iint \tr[ (u-u_\d)\otimes (u-u_\d) \cdot \n (u_\d \phi \chi_\e) ] \\
& + \iint \tr[ u_\d \otimes u_\d \cdot \n (u_\d \phi \chi_\e) ] \\
& = F_1+F_2+F_3,
\end{align}
where
$$
r_\d(u,u)(x) = \d^{-3} \int_{\R^3} \b(y \d^{-1}) (u(x-y) - u(x)) \otimes (u(x-y) - u(x)) dy.
$$
We will show that $F_1$ and $F_2$ already vanish in the limit of $\d \ra 0$ for a fixed $\e>0$. Let us observe the following estimate
\begin{equation}
 \|r_\d(u,u)\|_{9/4} \leq \d^{-3} \int_{\R^3} \b(y \d^{-1}) \|u(\cdot-y) - u(\cdot)\|^2_{9/2}dy
\overset{\text{def}}{=} R(t,\d).
\end{equation}
For $F_1$ we obtain
\begin{align*}
|F_1|& \leq \int_0^t R(t,\d) \|\n(u_\d\phi\chi_\e)\|_{9/5} dt \\ &\leq
\left(\int_0^t R^{3/2}(t,\d) dt \right)^{2/3} \left( \int_0^t \|\n(u_\d\phi\chi_\e)\|_{9/5}^3 dt \right)^{1/3}.
\end{align*}
Observe that
$$
R^{3/2}(t,\d) \leq \d^{-3} \int_{\R^3} \b(y \d^{-1}) \|u(\cdot-y) - u(\cdot)\|^3_{9/2}dy,
$$
and hence
$$
\int_0^t R^{3/2}(t,\d) dt \ra 0,
$$
as $\d \ra 0$, while in view of condition (ii) and \eqref{supp}, $\n(u_\d\phi\chi_\e) \in L^3L^{9/5}$ uniformly as $\d \ra 0$ for any fixed $\e>0$.
Thus, $F_1 \ra 0$. Similarly,
$$
|F_2| \leq \int_0^t \|u - u_\d\|_{9/2}^2  \|\n(u_\d\phi\chi_\e)\|_{9/5} dt \ra 0.
$$
As to $F_3$ we have
$$
F_3 = \frac{1}{2}\iint |u_\d|^2 u_\d \cdot \n(\phi \chi_\e) =
\frac{1}{2}\iint|u_\d|^2 \chi_\e u_\d \cdot \n \phi  +
\frac{1}{2}\iint |u_\d|^2 \phi\ u_\d \cdot \n \chi_\e .
$$
Clearly, the first integral on the right hand side converges to
$$
\frac{1}{2}\iint |u|^2 u\cdot \n \phi
$$
giving us the corresponding term in \eqref{mgee}. As to the second
integral we estimate using the  H\"{o}lder inequality and \eqref{der}
\begin{align*}
\left|\iint |u_\d|^2 \phi u_\d \cdot \n \chi_\e \right|& \leq \int_0^t \left(  \int_{|x-s_\e(t)| \leq 3\e}
|u_\d|^{9/2} dx \right)^{2/3} dt \\
& \leq \int_0^t \left(  \int_{|x-s_\e(t)| \leq 3\e+ \d}
|u|^{9/2} dx \right)^{2/3} dt \ra 0,
\end{align*}
as $\d,\e \ra 0$.

It remains to examine the pressure term $G$. We have
$$
G = \iint p_\d \chi_\e u_\d \cdot \n \phi + \iint p_\d \phi\ u_\d \cdot \n \chi_\e = G_1+G_2.
$$
Since $p \in L^{3/2}L^{9/4}$ and $u \in L^{3}L^{9/2}$ the
local $L^2$-pairing between $u$ and $p$ is continuous. So,
$$
G_1 \underset{\d\ra 0}{\ra} \iint p \chi_\e u \cdot \n \phi \underset{\e \ra 0}{\ra}
\iint p u\cdot \n \phi.
$$
As for $G_2$ we apply the following estimate
\begin{align*}
 |G_2| & \leq \|p_\d\|_{L^{3/2}L^{9/4}} \|\n \chi_\e\|_{L^\infty L^3}
\left[\int_0^t \left( \int_{|x-s_\e(t)|\leq 3\e} |u_\d|^{9/2} dx   \right)^{2/3} dt \right]^{1/3} \\
& \leq C \left[ \int_0^t\left( \int_{|x-s_\e(t)|\leq 3\e+\d} |u|^{9/2} dx   \right)^{2/3} dt \right]^{1/3}
\ra 0.
\end{align*}

This finishes the proof.

\section{Extentions}\label{S:ext}

First, we note that one can incorporate an external divergence-free force $f$ as long as
$\iint f \psi \ra \iint f u$ in the limit as $\d,\e \ra 0$. For this purpose $f \in H^{-1}$ appears to be sufficient.

\begin{figure}

\ifx\JPicScale\undefined\def\JPicScale{1}\fi
\psset{unit=\JPicScale mm}
\psset{linewidth=0.3,dotsep=1,hatchwidth=0.3,hatchsep=1.5,shadowsize=1,dimen=middle}
\psset{dotsize=0.7 2.5,dotscale=1 1,fillcolor=black}
\psset{arrowsize=1 2,arrowlength=1,arrowinset=0.25,tbarsize=0.7 5,bracketlength=0.15,rbracketlength=0.15}
\begin{pspicture}(0,0)(70,70)
\psline{->}(0,0)(70,0)
\psline{->}(0,0)(0,70)
\newrgbcolor{userHatchColour}{0.8 0.8 1}
\pscustom[linewidth=0.2,fillstyle=hlines,hatchwidth=0.1,hatchsep=1,hatchcolor=userHatchColour]{\psline(0,0)(33,0)
\psline(33,0)(25,25)
\psline(25,25)(0,50)
\psline(0,50)(0,0)
\psbezier[liftpen=2](0,0)(0,0)(0,0)(0,0)
}
\newrgbcolor{userHatchColour}{0.8 0.8 1}
\pscustom[linewidth=0.2,fillstyle=crosshatch,hatchwidth=0.1,hatchsep=1,hatchcolor=userHatchColour]{\psline(25,25)(22,33)
\psline(22,33)(0,60)
\psline(0,60)(0,50)
\psbezier[liftpen=2](0,50)(0,50)(0,50)(0,50)
}
\newrgbcolor{userHatchColour}{0.8 0.8 1}
\psline[linewidth=0.2,fillstyle=vlines,hatchwidth=0.1,hatchcolor=userHatchColour](50,0)(16,50)
\newrgbcolor{userHatchColour}{0.8 0.8 1}
\psline[linewidth=0.15,linestyle=dotted,dotsep=0.5,fillstyle=vlines,hatchwidth=0.1,hatchcolor=userHatchColour](22,33)(50,0)
\newrgbcolor{userHatchColour}{0.8 0.8 1}
\psline[linewidth=0.15,linestyle=dotted,dotsep=0.5,fillstyle=vlines,hatchwidth=0.1,hatchcolor=userHatchColour](22,33)(16,50)
\newrgbcolor{userHatchColour}{0.8 0.8 1}
\rput(73,0){\small{$\frac{1}{s}$}}
\newrgbcolor{userHatchColour}{0.8 0.8 1}
\rput(-4,69){\small{$\frac{1}{r}$}}
\newrgbcolor{userHatchColour}{0.8 0.8 1}
\rput{-45}(8,46){\tiny{new}}
\rput(14.4,12.3){Lions}
\psline[linewidth=0.15,linestyle=dotted,dotsep=0.5](50,0)(25,25)
\rput(50,-3){\tiny{$L^\infty L^2$}}
\rput(33,-3){\tiny{$L^\infty L^3$}}
\rput(16,52){\tiny{$L^2L^6$}}
\rput(22,24){\tiny{II}}
\rput(24,34){\tiny{I}}
\rput(-6,60){\tiny{$L^{5/3}L^\infty$}}
\rput(-6,50){\tiny{$L^2L^\infty$}}
\end{pspicture}

\bigskip

\caption{Here, $I = L^3L^{9/2}$, $II = L^4L^4$}\label{figa}

\end{figure}
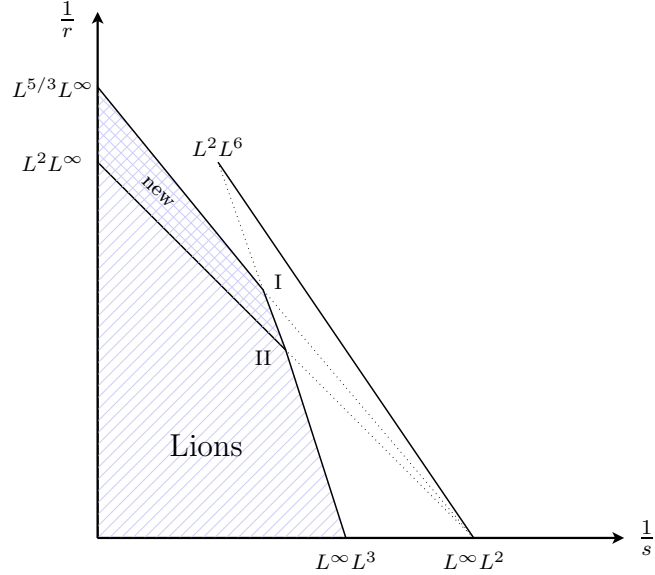

Second, by extrapolation from $L^3L^{9/2}$ along the line starting at $L^\infty L^2$ or
$L^2L^6$ or any other space in between we obtain the convex range of
$L^rL^s$-spaces determined by
\begin{equation}\label{extr}
 5/3r + 2/s \leq 1,\quad 1/r + 3/s \leq 1,
\end{equation}
for $r \geq 5/3$ and $s \geq 3$, which can be used in (i) to substitute
$L^3L^{9/2}$. Figure \ref{figa} graphically demonstrates the
region where this range is not covered by the previously known
results. The complementary condition on $\n u$ is given by
\begin{equation}\label{uext}
\n u \in L^{\frac{r}{r-2}}L^{\frac{s}{s-2}},
\end{equation}
for those $r,s>2$ that are in the range \eqref{extr}. For some $r$
and $s$, however, condition \eqref{uext} already implies (ii) by
interpolation with $\n u \in L^2L^2$ or it may be strong enough to
imply regularity via the Prodi-Serrin condition or \cite{Veiga}.
We leave details for the reader.

Treating the nonlinear terms $F$ as in \cite{ccfs} one can lower the
order of derivative in condition (ii) by
cost of increasing the integrability exponent. At extreme one
gets
\begin{equation}
u \in L^3B^{1/3}_{3,p}(([0,T]\times \R^3)\backslash
\mathrm{Graph(s)})_{\text{loc}},
\end{equation}
where "loc" means that for any $\phi \in D([0,T]\times
\R^3)\backslash \mathrm{Graph(s)})$, $\phi u \in
L^3B^{1/3}_{3,p}$.

We also notice that the argument does not make use of the global
estimates on $u$ and $p$. Thus if we are to pursue only the local
energy inequality \eqref{mgee} one can restate condition (i) in
the local sense with an additional assumption $p \in
L^{3/2}L^{9/4}_{\text{loc}}$. The letter does not seem to follow
directly from the corresponding conditions on $u$ without extra
smoothness assumptions on the initial condition (see \cite{SW}).

Lastly, we note that \thm{T:main} is valid on a smooth bounded
domain as well with the same requirement $p \in L^{3/2}L^{9/4}$.

\def\cprime{$'$}


\end{document}